\documentclass[reqno,12pt]{amsart}

\usepackage{amsmath}
\usepackage{amssymb}
\usepackage{amsthm}
\usepackage[english]{babel}
\newtheorem{theorem}{Theorem}

\newtheorem{lemma}{Lemma}

\newtheorem{corollary}{Corollary}

\begin{document}

\title {A new bound in the  Littlewood--Offord
problem}

\author[F.~G\"otze]{Friedrich G\"otze}
\author[A.Yu. Zaitsev]{Andrei Yu. Zaitsev}

\email{goetze@math.uni-bielefeld.de}
\address{Fakult\"at f\"ur Mathematik,\newline\indent
Universit\"at Bielefeld, Postfach 100131,\newline\indent D-33501 Bielefeld,
Germany\bigskip}
\email{zaitsev@pdmi.ras.ru}
\address{St.~Petersburg Department of Steklov Mathematical Institute
\newline\indent
Fontanka 27, St.~Petersburg 191023, Russia\newline\indent
and St.Petersburg State University, 7/9 Universitetskaya nab., St. Petersburg,
199034 Russia}

\begin{abstract}{The paper deals with studying a connection of the
Littlewood--Offord problem with estimating the concentration
functions of some symmetric infinitely divisible distributions. It is shown that the concentration function of a weighted sum of independent identically distributed random variables is estimated in terms of the concentration function of a symmetric infinitely divisible distribution whose spectral measure is concentrated on the set of plus-minus weights.
}\end{abstract}
\thanks{The authors were supported by SFB 1283/2 2021 –- 317210226 and by the RFBR-DFG grant
 20-51-12004.
The second author was supported by grant RFBR 19-01-00356.}

\maketitle

The aim of the present work is to provide a supplement to the  paper of Eliseeva and  Zaitsev~\cite{EZ14}. We study a connection of the
Littlewood--Offord problem with estimating the concentration
functions of some symmetric infinitely divisible distributions. In the proof, we repeat the
arguments of \cite{EZ14} adding at the last step an application of Jensen's inequality.

Let $X,X_1,\ldots,X_n$ be independent identically distributed
(i.i.d.) random variables.  The
 concentration function of a $\mathbf{R}^d$-dimensional random
vector $Y$ with distribution $F=\mathcal L(Y)$ is defined by the
equality
\begin{equation}
Q(F,\lambda)=\sup_{x\in\mathbf{R}^d}\mathbf{P}(Y\in x+ \lambda B),
\quad 0\le\lambda\leq\infty, \nonumber
\end{equation}
where $B=\{x\in\mathbf{R}^d:\|x\|\leq 1/2\}$. Of course, $Q(F,\infty)=1$.
Let $a=(a_1,\ldots,a_n)$, where
$a_k=(a_{k1},\ldots,a_{kd})\in \mathbf{R}^d$, $k=1,\ldots, n$. 
In this paper we
study the behavior of the concentration functions of the weighted
sums $S_a=\sum\limits_{k=1}^{n} X_k a_k$ with respect to the
properties of vectors~$a_k$. Interest in this subject
has increased considerably in connection with the study of
eigenvalues of random matrices (see, for instance, Friedland and
Sodin \cite{Fried:Sod:2007}, 
Rudelson and Vershynin \cite{Rud:Ver:2008}, \cite{Rud:Ver:2009},
Tao and Vu \cite{Tao:Vu:2009:Ann}, \cite{Tao:Vu:2009:Bull}, Nguyen and Vu~\cite{Nguyen:Vu:2011},
 Vershynin \cite{Ver:2011}, Tikhomirov \cite{T}, Livshyts,  Tikhomirov and Vershynin \cite{LTV}, Campos et al.
 \cite{CJMS}). For a detailed history of
the problem we refer to a review of Nguyen and Vu
\cite{Nguyen and Vu13}. The authors of the above articles (see
also Hal\'asz \cite{Hal:1977}) called this question the
Littlewood--Offord problem, since, for the first time, this
problem was considered in 1943 by Littlewood and Offord
\cite{Lit:Off:1943}  in connection with the study of random
polynomials. They considered a special case, where the
coefficients $a_k \in \mathbf{R}$ are one-dimensional, and $X$
takes values $\pm1$ with probabilities~$1/2$.

The recent achivements in estimating the probabilities of singularity of random matrices  \cite{T}, \cite{LTV}, 
 \cite{CJMS} were based on the Rudelson and Vershynin \cite{Rud:Ver:2008}, \cite{Rud:Ver:2009}, \cite{Ver:2011} method of {\it least common denominator.}
Note that the results of \cite{Fried:Sod:2007}, \cite{Rud:Ver:2009}, \cite{Ver:2011} (concerning the
Littlewood--Offord problem) were improved and refined in  \cite{Eliseeva and Zaitsev}, \cite{Eliseeva}, \cite{EGZ}.

Now we  introduce some notation. In the sequel, let $F_a$ denote
the distribution of the sum $S_a$, let $E_y$ be the probability
measure concentrated at a point $y$, and let $G$ be the
distribution of  symmetrized random variable
$\widetilde{X}=X_1-X_2$.
For $\delta\ge0$ we denote
\begin{equation}\label{8j88} p(\delta)=
G\big\{\{z:|z| > \delta\}\big\}.\end{equation}

The symbol $c$ will be used for absolute positive constants which
may be different even in the same formulas.

Writing $A\ll B$ means that $|A|\leq c B$. Also we will write
$A\asymp B$, if $A\ll B$ and $B\ll A$. We will write $A\ll _{d}
B$, if $|A|\leq c(d) B$, where $c(d)>0$ depends on $d$ only.
Similarly, $A\asymp_{d} B$, if $A\ll_{d} B$ and $B\ll_{d} A$. The
scalar product in $\mathbf{R}^d$ will be denoted $\left\langle \,
\cdot\,,\,\cdot\,\right\rangle$. Later $\lfloor x\rfloor$ is the
largest integer~$k$ such that~$k< x$.
 For~${x=(x_1,\dots,x_n )\in\mathbf R^n}$ we will use the norms
 $\|x\|^2= x_1^2+\dots +x_n^2$ and $|x|=
\max_j|x_j|$. We denote by $\widehat F(t)$, $t\in\mathbf R^d$, the
characteristic function of $d$-dimensional distributions~$F$.

Products and powers of measures will be understood in the
convolution sense. For infinitely
divisible distribution~$F$  and $\lambda\ge0$, we denote by $F^\lambda$ the infinitely divisible
distribution with characteristic function $\widehat F^\lambda(t)$.

The elementary properties of concentration functions are well
studied (see, for instance, \cite{Arak:Zaitsev:1986},
\cite{Hen:Teod:1980}, \cite{Petrov:1972}). It is known that
\begin{equation}\label{8jrt}
Q(F,\mu)\ll_d (1+ \lfloor \mu/\lambda \rfloor)^d\,Q(F,\lambda)
\end{equation}
 for any $\mu,\lambda>0$.
Hence,
\begin{equation}\label{8art}
Q(F,c\lambda)\asymp_{d}\,Q(F,\lambda).
\end{equation}

Let us formulate a generalization of the classical Ess\'een
inequality \cite{Ess:1966} to the multivariate case
(\cite{Ess:1968}, see also \cite{Hen:Teod:1980}):

\begin{lemma}\label{lm3} Let $\tau>0$ and let
 $F$ be a $d$-dimensional probability distribution. Then
\begin{equation}
Q(F, \tau)\ll_d \tau^d\int_{|t|\le1/\tau}|\widehat{F}(t)| \,dt.
\label{4s4d}
\end{equation}
\end{lemma}

In the general case $Q(F,\tau)$ cannot be estimated from below
by the right hand side of
 inequality~\eqref{4s4d}. However, if we assume additionally that the distribution $F$ is symmetric and
 its characterictic function is non-negative for all~$t\in\mathbf R$, then we have the lower bound:
\begin{equation} \label{1a}Q(F, \tau)\gg_d \tau^d\int\limits_{|t|\le1/\tau}{\widehat{F}(t)
\,dt},
\end{equation}
and, therefore,
\begin{equation} \label{1b}
Q(F, \tau)\asymp_d
\tau^d\int\limits_{|t|\le1/\tau}{\widehat{F}(t) \,dt},
\end{equation} (see \cite{Arak:1980} or \cite{Arak:Zaitsev:1986}, Lemma~1.5 of Chapter II for $d=1$).
In the multivariate case relations \eqref{1a} and~\eqref{1b} may be found in Zaitsev \cite{Zaitsev:1987}. Exactly the use of relation \eqref{1b} allows us to
simplify the arguments of Friedland and Sodin
 \cite{Fried:Sod:2007},  Rudelson and Vershynin~\cite{Rud:Ver:2009} and Vershynin \cite{Ver:2011}
 which were applied to  Littlewood--Offord problem
 (see  \cite{Eliseeva and Zaitsev}, \cite{Eliseeva}, \cite{EGZ}).
 \medskip

The main result of this paper is a general inequality which
reduces the estimation of concentration functions in the
Littlewood--Offord problem to the estimation of concentration
functions of some infinitely divisible distributions. This result
is formulated in Theorem~\ref{lm43}.

For $z\in \mathbf{R}$, introduce the distribution $H_z$ with the
characteristic function
\begin{equation} \label{11b}\widehat{H}_z(t)
=\exp\Big(-\frac{\,1\,}2\;\sum_{k=1}^{n}\big(1-\cos(\left\langle
\, t,a_k\right\rangle z)\big)\Big).
\end{equation}It depends on the vector~$a$.
It is clear that $H_z$ is a symmetric infinitely divisible
distribution.
  Therefore, its characteristic function is positive for all $t\in \mathbf{R}^d$.

Recall that $G=\mathcal L(X_1-X_2)$ and $F_a=\mathcal L(S_a)$, where $S_a=\sum\limits_{k=1}^{n} X_k a_k$.

\begin{theorem}\label{lm43}Let\/ $V$ be an arbitrary\/
one-dimensional Borel measure such that\/ $\lambda=V\{\mathbf
R\}>0$,\/ and $V\le G$, that is,\/ $V\{B\}\le G\{B\}$, for any
Borel set~$B$. Then, for any\/  $\tau>0$, we
have
\begin{equation}\label{cc23}
Q(F_a,\tau) \ll_d \int_{z\in\mathbf{R}}Q(H_{1}^{\lambda},\tau
|z|^{-1})\,W\{dz\},
\end{equation}where $W=\lambda^{-1}V$.\end{theorem}

\begin{corollary}\label{lm42}For any\/
$\varepsilon,\tau>0$, we have
\begin{equation}\label{11abc}
Q(F_a, \tau) \ll_d Q(H_1^{p(\tau/\varepsilon)}, \varepsilon),
\end{equation}
where $p(\,\cdot\,)$ is defined in \eqref{8j88}.
\end{corollary}

In order to verify Corollary \ref{lm42} we note that the
distribution $G=\mathcal{L}(\widetilde{X})$ may be represented as the mixture
$$
G=p_0G_0+p_1G_1,\quad\text{where}\quad p_j={\mathbf P}\bigl\{\widetilde{X} \in A_j\bigr\},\quad j=0,1,
$$
 $A_0=\{x\colon |x|\le\tau/\varepsilon\}$, $A_1=\{x\colon |x|>\tau/\varepsilon\}$,   $G_j$ are probability measures
defined for $p_j>0$ by the formula
$G_j\{B\}=G\{B\cap A_j\}/{p_j}$ , for any Borel set~$B$.
In fact, $G_j$ is the conditional distribution of $\widetilde X$
given that $\widetilde X\in A_j$. If $p_j=0$, then we can take as
$G_j$  an arbitrary  measure.

The conditions of  Theorem~\ref{lm43} are satisfied for $V=p_1G_1$. $\lambda=p_1=p(\tau/\varepsilon)$, $W=G_1$.

Inequalities  \eqref{8jrt} and \eqref{1b} imply that

\begin{eqnarray}Q(F_a,\tau) &\ll_d &\int_{z\in A_1}Q(H_{1}^{\lambda},\tau
|z|^{-1})\,W\{dz\}\nonumber \\
&\le&\sup_{z\ge\tau/\varepsilon}\;
Q\big(H_{1}^{p(\tau/\varepsilon)},\tau/z\big)=Q\big(H_{1}^{p(\tau/\varepsilon)},\varepsilon\big),
\end{eqnarray}
 proving~\eqref{11abc}.

Applying Theorem~\ref{lm43} with $V$ of the form
\begin{equation}\label{cc243}
V\{dz\}=\big(1+\lfloor\tau(\varepsilon
|z|)^{-1}\rfloor \big)^{-d}\,G\{dz\},
\end{equation}
and using inequality  \eqref{8jrt}, we obtain
\begin{corollary}\label{l429}
For any\/ $\varepsilon,\tau>0$, we have
\begin{equation}\label{c234}
Q(F_a,\tau) \ll_d
\lambda^{-1}\,Q(H_{1}^{\lambda},\varepsilon),
\end{equation}where
\begin{equation}\label{cc239}
\lambda=\lambda(G,\tau/\varepsilon)=V\{\mathbf
R\}=\int\limits_{z\in\mathbf{R}}\big(1+\lfloor\tau(\varepsilon
|z|)^{-1}\rfloor \big)^{-d}\,G\{dz\}.
\end{equation}
\end{corollary}

It is clear that $\lfloor\tau(\varepsilon
|z|)^{-1}\rfloor=0$ if $|z|>\tau/\varepsilon$. Therefore. $\lambda=\lambda(G,\tau/\varepsilon)\ge p(\tau/\varepsilon)$ and hence $Q(H_{1}^{\lambda},\varepsilon)\le Q(H_{1}^{p(\tau/\varepsilon)},\varepsilon)$. Thus, if $\lambda\gg_d 1$, then inequality~\eqref{c234} of Corollary~\ref{l429} is stronger than inequality~\eqref{11abc} of Corollary~\ref{lm42}.

The proof of Theorem~$\ref{lm43}$ is based on elementary
properties of concentration functions. We repeat the
arguments of \cite{EZ14} adding at the last step an application of Jensen's inequality. 
In \cite{EZ14}, inequality  \eqref{8jrt} was used instead. The main result of \cite{EZ14} does not imply Corollary~\ref{l429}.
Note that $H_1^{\lambda}$ is an infinitely divisible distribution
with the L\'evy
 spectral measure $M_\lambda=\frac{\,\lambda\,}4\;M^*$, where
 $M^*=\sum\limits_{k=1}^{n}\big(E_{a_k}+E_{-a_k}\big)$.
It is clear that the assertions of Theorem~$\ref{lm43}$ and
Corollaries~\ref{lm42} and~\ref{l429} may be treated as
statements about the measure~$M^*$.

Corollary \ref{lm42} was already proved earlier in
\cite{EZ14} and \cite{EGZ17}, see also \cite{Z15} for the case $\tau=0$. It was used essentially in \cite{EGZ17} and \cite{GZ18} to show that Arak's inequalities for concentration functions may be used for investigations
of  the Littlewood--Offord problem. Arak has shown that if the concentration function of infinitely divisible
distribution is relatively large, then the spectral measure of this distribution is concentrated in a neighborhood of a set with simple arithmetical structure. Together with Corollary \ref{lm42}, this means that large value of $Q(F_a,\tau)$ implies  a simple arithmetical structure of the set $\{\pm a_k\}_{k=1}^n$. This statement is similar to the so-called "inverse principle" in the Littlewood--Offord problem (see \cite{Tao:Vu:2009:Ann}, \cite{Nguyen:Vu:2011},  \cite{Nguyen and Vu13}).

Note that using the results of Arak
 \cite{Arak:1980},~\cite{Arak:1981} (see also \cite{Arak:Zaitsev:1986})
one could derive from Corollary~\ref{lm42} inequalities similar to
boumds for concentration functions in the Littlewood--Offord
problem, which were obtained in a  paper of Nguyen and
Vu~\cite{Nguyen:Vu:2011} (see also \cite{Nguyen and Vu13}). A
detailed discussion of this fact is presented in    \cite{EGZ17} and \cite{GZ18}. We noticed that Corollary~\ref{l429} may be stronger than  Corollary~\ref{lm42}. Therefore, the results of \cite{EGZ17} and \cite{GZ18} could be improved (in the sense of dependence of constants on the distribution of $X_1$) replacing an application of Corollary~\ref{lm42} by an application of Corollary~\ref{l429}.
The authors are going to devote  a separate publication to this topic.
\medskip

\noindent {\bf Proof of Theorem~\ref{lm43}.} Let us show that, for
arbitrary probability distribution~$W$ and $\lambda,T>0$,
\begin{multline}
\log\int_{|t|\le
T}\exp\Big(-\frac{\,1\,}2\;\sum_{k=1}^{n}\int_{z\in\mathbf{R}}\big(1-\cos(\left\langle
\, t,a_k\right\rangle z)\big)\,\lambda\,W\{dz\}\Big)\,dt\\ \le
\int_{z\in\mathbf{R}}\bigg(\log\int_{|t|\le
T}\exp\Big(-\frac{\,\lambda\,}2\;\sum_{k=1}^{n}\big(1-\cos(\left\langle
\, t,a_k\right\rangle z)\big)\Big)\,dt\bigg)\,W\{dz\}\\ =
\int_{z\in\mathbf{R}}\bigg(\log\int_{|t|\le
T}\widehat{H}_{z}^{\lambda}(t)\,dt\bigg)\,W\{dz\}.\label{dd11}
\end{multline}
It suffices to prove \eqref{dd11} for discrete  distributions~$W=
\sum_{j=1}^{\infty}p_j E_{z_j} $, where $0\le p_j\le1$,
$z_j\in\mathbf{R}$, $\sum_{j=1}^{\infty}p_j =1 $. Applying in this
case the generalized  H\"older inequality, we have
\begin{multline}
\int\limits_{|t|\le
T}\exp\Big(-\frac{\,1\,}2\;\sum\limits_{k=1}^{n}\int\limits_{z\in\mathbf{R}}\big(1-\cos(\left\langle
\, t,a_k\right\rangle z)\big)\,\lambda\,W\{dz\}\Big)\,dt\\ =
\int\limits_{|t|\le
T}\exp\Big(-\frac{\,\lambda\,}2\;\sum\limits_{j=1}^{\infty}p_j\sum\limits_{k=1}^{n}\big(1-\cos(\left\langle
\, t,a_k\right\rangle z_j)\big)\Big)\,dt\\ \le\prod\limits_{j=1}^{\infty}
\bigg(\int\limits_{|t|\le
T}\exp\Big(-\frac{\,\lambda\,}2\;\sum\limits_{k=1}^{n}\big(1-\cos(\left\langle
\, t,a_k\right\rangle z_j)\big)\Big)\,dt\bigg)^{p_j}.\label{d11}
\end{multline}
Taking logarithms of the left and right-hand sides
of~\eqref{d11}, we get~\eqref{dd11}. In general case we can
approximate~$W$ by discrete distributions in the
sense of weak convergence and pass to the limit.  Note also that
the integrals $\int_{|t|\le T}dt$ may be replaced in~\eqref{dd11} by
the integrals $\int\mu\{dt\}$ with an arbitrary Borel measure~$\mu$.

 Since for characteristic function
$\widehat{U}(t)$ of a random vector $Y$, we have
$$|\widehat{U}(t)|^2 = \mathbf{E}\exp(i\langle
\,t,\widetilde{Y}\rangle ) = \mathbf{E}\cos(\langle
\,t,\widetilde{Y}\rangle ),$$ where $\widetilde{Y}$ is the
corresponding symmetrized random vector, then
\begin{equation}\label{6}|\widehat{U}(t)| \leq
\exp\Big(-\cfrac{\,1\,}{2}\,\big(1-|\widehat{U}(t)|^2\big)\Big)  =
\exp\Big(-\cfrac{\,1\,}{2}\,\mathbf{E}\,\big(1-\cos(\langle
\,t,\widetilde{Y}\rangle )\big)\Big).
\end{equation}

According to Lemma \ref{lm3} and relations $V=\lambda\,W\le G$,
\eqref{dd11} and \eqref{6}, applying Jensen's inequality of the form $\exp(\mathbf{E}\,f(\xi))\le\mathbf{E}\,\exp(f(\xi))$ for any measurable function $f$ and any random varialble $\xi$, we have
\begin{eqnarray}
Q(F_a,\tau)&\ll_{d}& \tau^d\int_{\tau|t|\le 1}|\widehat{F_a}(t)|\,dt \nonumber\\
&\ll_{d}& \tau^d\int_{\tau|t|\le
1}\exp\Big(-\frac{\,1\,}{2}\,\sum_{k=1}^{n}\mathbf{E}\,\big(1-\cos(\left\langle
\,t,a_k\right\rangle \widetilde{X})\big)\Big)\,dt\nonumber\\
&=&\tau^d\int_{\tau|t|\le
1}\exp\Big(-\frac{\,1\,}2\;\sum_{k=1}^{n}\int_{z\in\mathbf{R}}\big(1-\cos(\left\langle
\, t,a_k\right\rangle z)\big)\,G\{dz\}\Big)\,dt\nonumber\\
&\le&\tau^d\int_{\tau|t|\le
1}\exp\Big(-\frac{\,1\,}2\;\sum_{k=1}^{n}\int_{z\in\mathbf{R}}\big(1-\cos(\left\langle
\, t,a_k\right\rangle z)\big)\,\lambda\,W\{dz\}\Big)\,dt\nonumber\\
&\le&\exp\bigg(
\int_{z\in\mathbf{R}}\log\bigg(\tau^d\int_{\tau|t|\le
1}\widehat{H}_{z}^{\lambda}(t)\,dt\bigg)\,W\{dz\}\bigg)\nonumber\\
&\le& \int_{z\in\mathbf{R}}\bigg(\tau^d\int_{\tau|t|\le
1}\widehat{H}_{z}^{\lambda}(t)\,dt\bigg)\,W\{dz\}.\label{cc11}
\end{eqnarray}

Using  \eqref{1b}, we have
\begin{eqnarray}
 \tau^d\int_{\tau|t|\le 1}
\widehat{H}_{z}^{\lambda}(t) \,dt &\asymp_{d}&
Q(H_{z}^{\lambda},\tau) = Q\big(H_{1}^{\lambda},
\tau{|z|}^{-1}\big)
.\label{cc22}
\end{eqnarray}
 Substituting this formula into
\eqref{cc11}, we obtain~\eqref{cc23}. In \eqref{cc22}, we have used that $H_z^{\lambda}=\mathcal L(z\eta) $, where $\eta$ is a random vector with $\mathcal L(\eta)=H_1^{\lambda} $. $\square$\medskip

{\bf Acknowledgment.} We are grateful to the reviewers for valuable remarks.

\end{document}